\documentclass{amsart}
\usepackage{amsfonts}
\usepackage{amsmath,amssymb}
\usepackage{amsthm}
\usepackage{amscd}
\usepackage{graphics}
\usepackage{graphicx}

\theoremstyle{remark}{
\newtheorem{Def}{{\rm Definition}}

\newtheorem{Rem}{{\rm Remark}}

}
\theoremstyle{plain}
{

\newtheorem{Prop}{Proposition}

\newtheorem{MainThm}{Main Theorem}

}

\begin{document}
\title[Shapes of regions formed by cylinders of circles of fixed radii]{Regions surrounded by cylinders of circles of fixed radii and exposition of their shapes by natural graphs}
\author{Naoki kitazawa}

\address{Institute of Mathematics for Industry, Kyushu University, 744 Motooka, Nishi-ku Fukuoka 819-0395, Japan\\
 TEL (Office): +81-92-802-4402 \\
 FAX (Office): +81-92-802-4405 \\
}
\email{naokikitazawa.formath@gmail.com}
\urladdr{https://naokikitazawa.github.io/NaokiKitazawa.html}
\maketitle
\begin{abstract}
We investigate regions formed by cylinders of circles of fixed radii. We investigate graphs obtained by collapsing each level set of the functions represented by the natural projections of them to the $1$-dimensional line.
Some specific trees obtained in simple ways from so-called balanced trees are shown to be realized as such graphs. Related studies on regions in the Euclidean plane surrounded by real algebraic curves are presented by several researchers. One of pioneering studies is presented by Bodin, Popescu-Pampu and Sorea in 2022--3 as an elementary and surprisingly new study. The author has been interested in related studies and also in constructing natural and explicit real algebraic maps onto such regions, generalizing the canonical projections of the unit spheres. Such studies in real algebraic geometry, different from theory of existence in the last century, mainly studied by Nash and Tognoli, are remarked.
\end{abstract}
%【REVISE】 combinatoric ～ is → combinatorial object. It is .
%【REVISE】  such that a point is a vertex if and only if the corresponding connected component of the level set contains some singular points → whose vertex set is the set of all points containing some singular points in the corresponding connected component of the level set .
%【REVISE】 We delete "extending the result before".
\section{Introduction.}
\label{sec:1}
Regions in the $k$-dimensional Euclidean space ({\it real affine space}) ${\mathbb{R}}^k$ surrounded by hypersurfaces are fundamental geometric and combinatorial objects. We focus on their shapes. We understand the shapes by graphs they collapse respecting the canonical projection to a $1$-dimensional affine line. The graph is the {\it Poincar\'e-Reeb graph} of the region. More precisely, the space obtained by contracting components of level sets to single points of the restriction of the projection. 

Related studies seem to be fundamental, natural, and surprisingly developing recently. One of pioneering studies is \cite{bodinpopescupampusorea}.
Their, regions in the plane surrounded by real algebraic curves and the graphs explained before are studied. A graph generically embedded into the plane is realized as the graph of a suitably found region by approximation of curves by the zero sets of polynomials.
Motivated by this, the author has contributed to development of related studies (essentially by himself). More precisely, the author has been interested in construction and examples in real algebraic geometry, different from existence and approximation mainly pioneered by Nash and Tognoli. For related existence theory, see \cite{kollar,nash, tognoli}.

\cite{kitazawa3} is a related pioneering study on the construction. This is motivated by \cite{bodinpopescupampusorea}. The author hit on an idea that this is regarded as the image of a nice smooth map of a certain class generalizing the canonical projections of the unit spheres, the class of {\it special generic} maps. For special generic maps, see also \cite{saeki1}. Related to this, the author has been studying higher dimensional versions of Morse functions and applications of differential topology of manifolds, or so-called global singularity theory and applications to differential topology of manifolds. For related studies on special generic maps, see \cite{kitazawa15} and another case, {\it round fold} maps, generalizing the canonical projections of the unit spheres by respecting the images of the set of all singular points of the maps, see \cite{kitazawa0.1, kitazawa0.2, kitazawa0.3}, where round fold maps are another topic in global singularity theory with differential topological applications. The author presents explicit real algebraic construction of a special generic map in this real algebraic situation. See also \cite{kitazawa8}, as a kind of additional remarks.
We add exposition on \cite{kitazawa3}. The {\it Reeb graph} of a smooth function is the space represented as a quotient space of the manifold of the domain. It is, same as the {\it Poincar\'e-Reeb graph}, obtained by considering all components of all level sets and contracting each component to a single point. Note that this notion has already appeared in \cite{reeb} and is classical and much older than the notion of {\it Poincar\'e-Reeb graph}. The notion of Poincar\'e-Reeb graph appeared recently and \cite{bodinpopescupampusorea, sorea1, sorea2} are related pioneering studies. {\it Reeb graphs} have been fundamental and strong tools in theory of Morse functions and related differential topology and recently, in several scenes where mathematics is applied, such as visualization. In \cite{kitazawa3}, the author has also investigated the Reeb graph of the function represented as the composition of the map with the canonical projection to the $1$-dimensional line. This is also motivated by reconstructing nice smooth functions on closed manifolds with prescribed Reeb graphs and level sets, mainly \cite{sharko}, followed by \cite{masumotosaeki, michalak}, and later, by \cite{kitazawa1, kitazawa2}.
For related study, see \cite{kitazawa4, kitazawa5, kitazawa6, kitazawa7, kitazawa12, kitazawa13} for example.

Related to the construction of real algebraic maps (of non-positive codimensions) by the author, the author has been interested in the regions surrounded by circles and curves in the plane and hypersurfaces in a general real affine space. This is essentially started by \cite{kitazawa11} with \cite{kitazawa9, kitazawa10}, where only the cases of curves in the plane are studied. Recently \cite{kitazawa14} is closely related. There regions surrounded by real algebraic curves, mainly circles of fixed radii, in the real affine plane, and ({\it refined}) {\it Poincar\'e-Reeb graphs} of them are discussed and investigated.

In the present paper, we first define regions surrounded by hypersurfaces in a general real affine space and its {\it Poincar\"e-Reeb graph}. This is not so difficult. Our main theorem is as follows. We explain explicit classes of graphs rigorously, later.
 
\setcounter{MainThm}{-1}
\begin{MainThm}
\label{mthm:0}
Let $G$ be a finite graph as follows{\rm :} prepare a balanced tree whose root is $v_0$, whose depth is $d \geq 1$, and the number of children of each vertex of a same depth $0 \leq i \leq d-1$ is always $n_{c,i+1} \geq 2$. Add a new edge connecting $v_0$ and another new vertex $v_1$ to this balanced tree. In this situation, $G$ is realized as the Poincar\'e-Reeb graph of a region surrounded by finitely many copies of a cylinder of a circle of some suitably chosen fixed radius in some real affine space. 
\end{MainThm}
This paper is organized as follows. The next section is for preliminaries and we introduce or review fundamental notions and notation on topological spaces, smooth or {\it real algebraic} manifolds (curves and hypersurfaces), maps between the manifolds, and (di)graphs, more rigorously. 
We also explain {\it Reeb graphs}, {\it Reeb digraphs} and {\it Reeb V-digraphs}. We also introduce the notion of {\it RA-region}. An {\it RA-region} is a region surrounded by finitely many real algebraic hypersurfaces intersecting in the transversal way in a real affine space. Such regions of dimension greater than $2$ appear in several studies such as \cite{kitazawa4, kitazawa5, kitazawa6, kitazawa7, kitazawa12}, by the author. We also define the {\it Poincar\'e-Reeb V-digraph} of it and this is done first in the present paper.
In the third section, we present our main result again in a more rigorous and refined way and prove it. In the fourth section, we also explain related studies such as reconstruction of real algebraic maps onto a region by the author and some previous result on regions surrounded by real algebraic hypersurfaces and their Poincar\'e-Reeb graphs by the author.\\
\ \\
\noindent {\bf Conflict of interest.} \\
The author is also a researcher at Osaka Central
Advanced Mathematical Institute (OCAMI researcher), supported by MEXT Promotion of Distinctive Joint Research Center Program JPMXP0723833165. He is not employed there. This is for our studies
and our study also thanks this. \\
\ \\
{\bf Data availability.} \\
No other data is associated to the paper.
\section{Preliminaries.}
\subsection{Fundamental terminologies and notation.}
\subsubsection{Topological spaces, and smooth manifolds and maps.}
Let $\emptyset$ denote the empty set. 

For a topological space $X$ and its subspace $Y$, let ${\overline{Y}}^X$ denote the closure of $Y$ in $X$. For a topological space $X$ having the structure of a cell complex the dimensions of cell of which are bounded, the dimension $\dim X$ is defined uniquely as a non-negative integer and this yields a topological invariant. 
A topological manifold is well-known to have the structure of a CW complex. A smooth manifold has the structure of a polyhedron and the structure of a certain polyhedron is well-known as a {\it PL} manifold. A topological space having the structure of a polyhedron of dimension at most $2$ has the structure of a polyhedron uniquely. This is true for topological manifolds of dimension at most $3$, according to \cite{moise}. Hereafter, ${\rm Int}\ X$ is the interior of a (topological) manifold $X$ and $\partial X:=X-{\rm Int}\ X$ is its boundary.
 
The $k$-dimensional Euclidean space ${\mathbb{R}}^k$ is a simplest smooth manifold and the Riemannian manifold equipped with the so-called standard Euclidean metric, where $\mathbb{R}:={\mathbb{R}}^1$.
For each point $x \in {\mathbb{R}}^k$, $||x|| \geq 0$ is the distance between $x$ and the origin $0 \in {\mathbb{R}}^k$ under the metric. The $k$-dimensional unit sphere $S^k:=\{x \in {\mathbb{R}}^{k+1} \mid ||x||=1\}$ is a $k$-dimensional smooth compact submanifold of ${\mathbb{R}}^{k+1}$ with no boundary, which is connected for $k \geq 1$ and a discrete two-point set for $k=0$. It is also the zero set of the real polynomial $||x||^2-1={\Sigma}_{j=1}^{k+1} {x_j}^2-1$ with $x:=(x_1,\cdots,x_{k+1})$. The $k$-dimensional unit disk $D^k:=\{x \in {\mathbb{R}}^{k} \mid ||x|| \leq 1\}$ is a $k$-dimensional smooth compact and connected submanifold of ${\mathbb{R}}^{k}$ and we easily have $\partial D^k=S^{k-1}$.

For a $1$-dimensional smooth manifold, we also use "{\it curve}". For a submanifold $Y$ of dimension $m-1$ in an $m$-dimensional manifold $X \supset Y$, we also use "{\it hypersurface}". 

For a differentiable manifold $X$, the tangent vector space $T_x X$ at $x \in X$ is also defined as a real vector space of dimension $\dim X$. Let $c:X \rightarrow Y$ be a differentiable map from a differentiable manifold $X$ into another manifold $Y$. The differential ${dc}_x:T_x X \rightarrow T_{c(x)} Y$ of $c$ at $x \in X$ is a linear map. If the rank of ${dc}_x$ is smaller than the minimum between $\dim X$ and $\dim Y$, then $x$ is said to be a {\it singular point} of $c$, with $c(x)$ being a {\it singular value} of $c$. Let $S(c)$ denote the set of all singular points (the {\it singular set} of $c$). We also use "{\it critical}" instead of "singular" in the case of a real-valued function $c:X \rightarrow \mathbb{R}$. Hereafter, differentiable maps are smooth maps or maps of the class $C^{\infty}$. A {\it diffeomorphism} means a homeomorphism with no critical point and we can define the notion of {\it diffeomorphic} manifolds or equivalently, a notion that a manifold is diffeomorphic to another manifold. 
The canonical projection of the Euclidean space ${\mathbb{R}}^k$ into ${\mathbb{R}}^{k_1}$ is ${\pi}_{k,k_1}:{\mathbb{R}}^{k} \rightarrow {\mathbb{R}}^{k_1}$ with ${\pi}_{k,k_1}(x)=x_1$ where $x=(x_1,x_2) \in {\mathbb{R}}^{k_1} \times {\mathbb{R}}^{k_2}={\mathbb{R}}^k$ with $k_1, k_2>0$ and $k=k_1+k_2$. The canonical projection of the unit sphere $S^{k-1}$ is its restriction. 
\subsubsection{Real algebraic objects.}
A {\it real algebraic manifold} means a union of connected components of the zero set of a real polynomial map and a set which is also {\it non-singular}{\rm : }it is defined by the implicit function theorem for the real polynomial map. The space ${\mathbb{R}}^k$, which is also called the {\it $k$-dimensional real affine space}, and $S^{k-1} \subset {\mathbb{R}}^k$ are of simplest examples. {\it Real algebraic} maps are the compositions of the embeddings into the real affine spaces (defined canonically) with the canonical projections ${\pi}_{k,k_1}$.

A {\it circle} means $S^1$ or a circle of a fixed radius $r>0$ centered at $(x_{0,1},x_{0,2}) \in {\mathbb{R}}^2$  in ${\mathbb{R}}^2$: in other words, the zero set of the form $\{(x_1,x_2) \mid {(x_1-x_{0,1})}^2+{(x_2-x_{0,2})}^2-r=0\}$ with $(x_{0,1},x_{0,2}) \in {\mathbb{R}}^2$ and $r>0$.

A {\it cylinder} of a real algebraic manifold $C \subset {\mathbb{R}}^{k_1}$ in ${\mathbb{R}}^k={\mathbb{R}}^{k_1} \times {\mathbb{R}}^{k_2}$ is a real algebraic manifold of the form $C \times {\mathbb{R}}^{k_2} \subset {\mathbb{R}}^k$. 
The hyperplane $\{0\} \times {\mathbb{R}}^{k_2}={\mathbb{R}}^{k_2+1}={\mathbb{R}}^k$ with $k_1=1$ is a cylinder of a single point in $\mathbb{R}$, which is the zero set of a real polynomial of degree $1$. A subset obtained by a finite iteration of transformations each of which is a parallel transformation, a rotation around a point, or a reflection around the hyperplane $\{0\} \times {\mathbb{R}}^{k_2} \subset {\mathbb{R}}^{k_2+1}={\mathbb{R}}^{k}$, is also a real algebraic manifold of dimension $\dim C+k_2$. If $C$ is of a class $\mathcal{C}_{\rm R.A}$ of real algebraic manifolds in ${\mathbb{R}}^{k_1}$, then these real algebraic manifolds are called cylinders of {\it the class $\mathcal{C}_{\rm R.A}$}. For example, if $\mathcal{C}_{\rm R.A}$ is the class of single points in $\mathbb{R}$, then the resulting real algebraic manifolds are real algebraic hypersurfaces in ${\mathbb{R}}^{k_2+1}$ and called cylinders of single points. In such a case, {\it real algebraic hyperplanes}, the zero sets of real polynomial functions of degree $1$, are obtained. As another case, if $\mathcal{C}_{\rm R.A}$ is the class of circles, then the resulting real algebraic hypersurfaces are called the cylinders of circles, and the zero sets of real polynomial functions of degree $2$ of a certain simple form.

We mainly consider circles and cylinders of a circle in our main result. By considering the disk $D_C \subset {\mathbb{R}}^2$ with $\partial D_C=C$ being a circle, we can define a {\it cylinder of a disk}, whose boundary is a cylinder of a circle (C). 
The complementary set is a cylinder of {\it the complementary set of a disk}.
We consider these notions in Main Theorems \ref{mthm:1} and \ref{mthm:2}.
  
\subsubsection{Graphs, digraphs, V-digraphs and Reeb {\rm (}V-di{\rm )}graphs.}
A graph $G$ means a $1$-dimensional finite and connected cell (CW) complex such that the closure of each $1$-cell is homeomorphic to $D^1$.
A {\it vertex} (an {\it edge}) of a graph is a 0-cell (resp. 1-cell) and the set of all vertices (resp. edges) of the graph is the {\it vertex set} (resp. {\it edge set}) of the graph.

The notion of {\it isomorphic} graphs, or equivalently the notion that a graph is isomorphic to the other graph, is defined by the existence of a homeomorphism mapping the vertex set of one graph to that of the other graph.

The {\it degree} $\deg(v)$ of a vertex $v$ of $G$ is the number of edges $e$ incident to $v$, or $e \ni v$. 

A {\it non-trivial} graph is a graph which is not isomorphic to a graph with exactly one vertex.
A {\it tree} is a non-trivial graph whose 1st Betti number is $0$.

A {\it digraph} $\overrightarrow{G}$ is a graph with each edge of it being oriented. We can define these notions for graphs in the same way.
The {\it indegree} ({\it outdegree}) ${\deg}_{\rm in}(v)$ (resp. ${\deg}_{\rm out}(v)$) of a vertex $v$ of $\overrightarrow{G}$ is the number of edges $e$ entering (resp. departing from) $v$.
We call a vertex $v \in \overrightarrow{G}$ {\it source} ({\it sink}) if ${\deg}_{\rm in}(v)=0$ (resp. ${\deg}_{\rm out}(v)=0$).
\begin{Prop}
\label{prop:1}
A tree $G$ has the structure of a digraph by the following. 
\begin{itemize}
\item We can choose a vertex $v_0 \in G$.
\item We can define the structure of a digraph $\overrightarrow{G}$ on $G$ in such a way that at exactly
$l>0$ vertices $v_{1,j}$, labeled by integers $1 \leq j \leq l$, ${\deg}_{\rm out}(v_{1,j})=0$ and ${\deg}_{\rm in}(v_{1,j})=1$,
except $l>0$ vertices $v_{1,j}$, at each vertex $v$, ${\deg}_{\rm out}(v) \geq {\deg}_{\rm in}(v)=1$ and that especially, $v_0$ is the unique source of $\overrightarrow{G}$.
\end{itemize}
\end{Prop}
\begin{Def}
\label{def:1}
In the situation of Proposition \ref{prop:1}, the triplet $(\overrightarrow{G},v_0,\{v_{1,j}\}_{j=1}^{l})$ is called a {\it rooted tree}
where $\{v_{1,j}\}_{j=1}^{l}$ is the family of all vertices $v_{1,j}$ with ${\deg}_{\rm in}(v_{1,j})=1>{\deg}_{\rm out}(v_{1,j})=0$ and all sinks of the digraph. In addition, the vertex $v_0$ is said to be the {\it root} of the rooted tree and each vertex $v_{1.j}$ is said to be a {\it leaf} of the rooted tree.
\begin{enumerate}
\item \label{def:1.1} In this situation, for vertices $v$ except vertices of $\{v_{1,j}\}_{j=1}^{l}$, ${\deg}_{\rm out}(v)$ is also called the {\it number of children} of $v$ (in the rooted tree).   
\item \label{def:1.2} For each vertex $v$, we can have a family of all edges of $\overrightarrow{G}$ to reach $v$ in the straight way from the root $v_0$ uniquely. The {\it depth} of $v$ (in the rooted tree) is defined as the number of these edges and a positive integer {\rm :} the depth of $v_0$ is $0$ and for remaining vertices $v \neq v_0$ in the rooted tree, the depth of $v$ is positive.
\item \label{def:1.3} Suppose that the depth of $v_{1,j}$ is always $d>0$ for each vertex $v_{1,j}$. Suppose also that the number of children of a vertex $v$ of depth $1 \leq i<d$ is always $n_{c,i+1} \geq 1$ and that the number of children of the root $v_0$ is $n_{c,1}$, then the rooted tree is called a {\it balanced tree of type $(\{n_{c,i}\}_{i=1}^{d},d)$}{\rm :} $d$ is defined to be the {\it depth} of the rooted tree. As a specific case, the case $n_{c,i}=1$ for $1 \leq i \leq d$ is for a so-called {\it path digraph on exactly $d+1$ vertices}. 
\end{enumerate}   
\end{Def}

A {\it V-digraph} $\overrightarrow{G_g}$ is a pair of a graph $G$ with a continuous function $g:G \rightarrow \mathbb{R}$ on the underlying graph $G$ which is on each edge injective and which induces the structure of a digraph $\overrightarrow{G_g}$ based on the orders of the values of $g$ at vertices of $G$ canonically.

For a smooth function $c:X \rightarrow \mathbb{R}$ on a closed manifold, we can define the following equivalence relation ${\sim}_c$ on $X$: $p_1 {\sim}_c p_2$ if and only if they are in a same component of some level set $c^{-1}(q)$. We have the quotient space $R_c$ with the quotient map $q_c:X \rightarrow R_c$ and the unique continuous function $\bar{c}:R_c \rightarrow \mathbb{R}$ enjoying the relation $c=\bar{c} \circ q_c$. 
In the case $c(S(c))$ is finite, then $R_c$ is a graph whose vertex set consists of all points corresponding to components of level sets of $c$ containing some points of the singular set $S(c)$ of $c$. 
This is the {\it Reeb graph} of $c$. We can induce the structure of a V-digraph $\overrightarrow{{R_c}_{\bar{c}}}$. This is the {\it Reeb {\rm (}V-{\rm )}digraph} of $c$. These notions do not appear directly in our main result. However, they are important in our important notion, the {\it Poincar\'e-Reeb V-digraph} of a region surrounded by real algebraic hypersurfaces (an {\it RA-region}). The notion of Reeb V-digraph appears again in the fourth section.
 
\subsection{An RA-region and its Poincar\'e-Reeb V-digraph.}
We define the notion of {\it RA-region}. This is a region surrounded by real algebraic hypersurfaces in the real affine space ${\mathbb{R}}^k$ intersecting in the transversal way in the closure of the region.
\begin{Def}
\label{def:2}
Let $D \subset {\mathbb{R}}^k$ be a non-empty, connected and open set of ${\mathbb{R}}^k$. Let $\{S_j\}_{j=1}^{l}$ be a family of non-singular real algebraic hypersurfaces in ${\mathbb{R}}^k$ such that $S_j$ is a union of some components of the zero set of a real polynomial $f_j$ and that $S_j  \bigcap {\overline{D}}^{{\mathbb{R}}^k}  \neq \emptyset$ for all integers $1 \leq j \leq l$. We also assume the following.
\begin{enumerate}
\item \label{def:2.1} There exists a small open neighborhood $U_D$ of ${\overline{D}}^{{\mathbb{R}}^k}$ such that $D=U_D \bigcap {\bigcap}_{j=1}^{l} \{x \mid f_j(x)> 0\}$ and $U_D \bigcap S_j=U_D \bigcap \{f_j(x)=0\}$ hold.
\item \label{def:2.2} {\rm (}On transversality{\rm )} For any increasing subsequence $\{i_j\}_{j=1}^{l^{\prime}}$ of the sequence $\{j\}_{j=1}^l$ of positive integers and any point  $p \in {\bigcap}_{j=1}^{l^{\prime}} S_{i_j} \bigcap {\overline{D}}^{{\mathbb{R}}^k}$, the relation $\dim {\bigcap}_{j=1}^{l^{\prime}} T_p S_{i_j}=k-l^{\prime}$ holds.
\end{enumerate}
The pair $(D,\{S_j\}_{j=1}^l)$ is called a {\it real algebraic region} or an {\it RA-region} of ${\mathbb{R}}^k$.  
\end{Def}
For the restriction ${\pi}_{k,1} {\mid}_{{\overline{D}}^{{\mathbb{R}}^k}}$, we can define the quotient space ${PR}_{(D,\{S_j\}_{j=1}^l)}$ and the quotient map $q_{(D,\{S_j\}_{j=1}^l)}:{\overline{D}}^{{\mathbb{R}}^k} \rightarrow {PR}_{(D,\{S_j\}_{j=1}^l)}$ and can have the unique continuous function $\overline{{\pi}_{(D,\{S_j\}_{j=1}^l)}}:{PR}_{(D,\{S_j\}_{j=1}^l)} \rightarrow \mathbb{R}$ with ${\pi}_{k,1} {\mid}_{{\overline{D}}^{{\mathbb{R}}^k}}=\overline{{\pi}_{(D,\{S_j\}_{j=1}^l)}} \circ q_{(D,\{S_j\}_{j=1}^l)}$ as we do in defining the Reeb V-digraph of a smooth function on a closed manifold.

We define the notion of {\it singular point} for an RA-region $\overline{{\pi}_{(D,\{S_j\}_{j=1}^l)}}$. A point $p \in {\overline{D}}^{{\mathbb{R}}^k}$ is {\it singular} if and only if either of the following two holds.
\begin{itemize}
\item The point $p$ is in exactly $l^{\prime}<k$ hypersurfaces $S_{i_j}$. We define the set $S_{i_j,1 \leq j \leq l^{\prime}}$ of all points of ${\overline{D}}^{{\mathbb{R}}^k}$ contained in these $l^{\prime}$ hypersurfaces $S_{i_j}$ and not contained in any other hypersurface $S_j$, which is seen to be a ($k-l^{\prime}$)-dimensional manifold with no boundary. The point $p$ is also a critical point of the restriction ${\pi}_{k,1} {\mid}_{S_{i_j,1 \leq j \leq l^{\prime}}}$.
\item  The point $p$ satisfies $p \in {\bigcap}_{j=1}^{l^{\prime}} S_{i_j}$ with $l^{\prime}=k$.
\end{itemize}

This is a generalization of arguments presented in \cite{kitazawa9, kitazawa10} for example. There the case $k=2$ is discussed.
Hereafter, ${\overline{D}}^{{\mathbb{R}}^k}$ is assumed to be compact.
Around a point $p$ of ${PR}_{(D,\{S_j\}_{j=1}^l)}$ which is corresponding to a component with no singular point of the RA-region, the quotient map is regarded as the projection of a product bundle over a small neighborhood of $p$ in ${PR}_{(D,\{S_j\}_{j=1}^l)}$. This is due to the so-called relative version of Ehresmann's fibration theorem or more generally, so-called stratified Morse theory. 
For stratified Morse theory, see \cite{massey} and see also \cite{hamm}.
By a kind of finiteness due to the real algebraic situation, we can define the structure of a graph on ${PR}_{(D,\{S_j\}_{j=1}^l)}$ whose vertex set consists of all points corresponding to components containing some singular points of $(D,\{S_j\}_{j=1}^l)$. This is the {\it Poincare-Reeb graph} of $(D,\{S_j\}_{j=1}^l)$. We also have the {\it Poincar\'e-Reeb {\rm (}V-{\rm )}digraph} $\overrightarrow{{{PR}_{(D,\{S_j\}_{j=1}^l)}}_{\overline{{\pi}_{(D,\{S_j\}_{j=1}^l)}}}}$ of $(D,\{S_j\}_{j=1}^l)$.

We can also understand these arguments from general theorems \cite[Theorem 3.1]{saeki2} and \cite[Theorem 2.1]{saeki3} and related arguments.

\section{Main Theorems and their proof.}
Main Theorem \ref{mthm:0} is revisited in a refined way.
\begin{MainThm}
\label{mthm:1}
Let $d>0$ be an integer. 
Let $\{t_j\}_{j=1}^{d+2}$ be an increasing sequence of real numbers{\rm :} $t_j<t_{j+1}$ for $1 \leq j \leq d+1$.
Let $\{n_{c,i}\}_{i=1}^{d}$ be a sequence of integers greater than $1$.
Let
$\overrightarrow{G}$ be a digraph obtained in the following way.
\begin{itemize}
\item Choose a balanced tree of type $(\{n_{c,i}\}_{i=1}^{d},d)$.
\item Add a new vertex and a new edge entering the root of the balanced tree above.
\end{itemize}
Then an RA-region $(D,\{S_j\}_{j=1}^{l}):=(D,\{S_j\}_{j=1}^{{\Sigma}_{i=1}^{d} (n_{c,i}+1)+2})$ surrounded by cylinders $S_j$ of a circle $C_R:=\{x=(x_1,x_2) \mid ||x||=R\}$ in ${\mathbb{R}}^2$ whose radius is any real number $R$ greater than a sufficiently large positive number $R_0$ can be chosen in such a way that the Poincar\'e-Reeb digraph $\overrightarrow{{{PR}_{(D,\{S_j\}_{j=1}^l)}}_{\overline{{\pi}_{(D,\{S_j\}_{j=1}^l)}}}}$ is isomorphic to $\overrightarrow{G}$ and that the image of its vertex set by the map $\overline{{\pi}_{(D,\{S_j\}_{j=1}^l)}}:{PR}_{(D,\{S_j\}_{j=1}^l)} \rightarrow \mathbb{R}$ is $\{t_j\}_{j=1}^{d+2}$.
\end{MainThm}
\begin{proof}
A main ingredient is to choose cylinders of a circle $C_R:=\{x=(x_1,x_2) \mid ||x||=R\}$ in ${\mathbb{R}}^2$ whose radius is $R>R_0$ (for the sufficiently large real number $R_0>0$) in a suitable way. Let $D_R \subset {\mathbb{R}}^2$ be a disk with $\partial D_R=C_R$. We have $D_R:=\{x=(x_1,x_2) \mid ||x|| \leq R\}$ of course.
We can choose the following. 
\begin{itemize}
\item The first two cylinders $S_1$ and $S_2$ of $C_R$ are as presented in Figure \ref{fig:1}. 
We can have a cylinder $D_{S_i}$ of the disk $D_R$ with $\partial D_{S_i}=S_i$ ($i=1,2$). 
The region represented as the intersection $D_{S_1} \bigcap D_{S_2}$ is also depicted.
\begin{figure}
	\includegraphics[width=80mm,height=80mm]{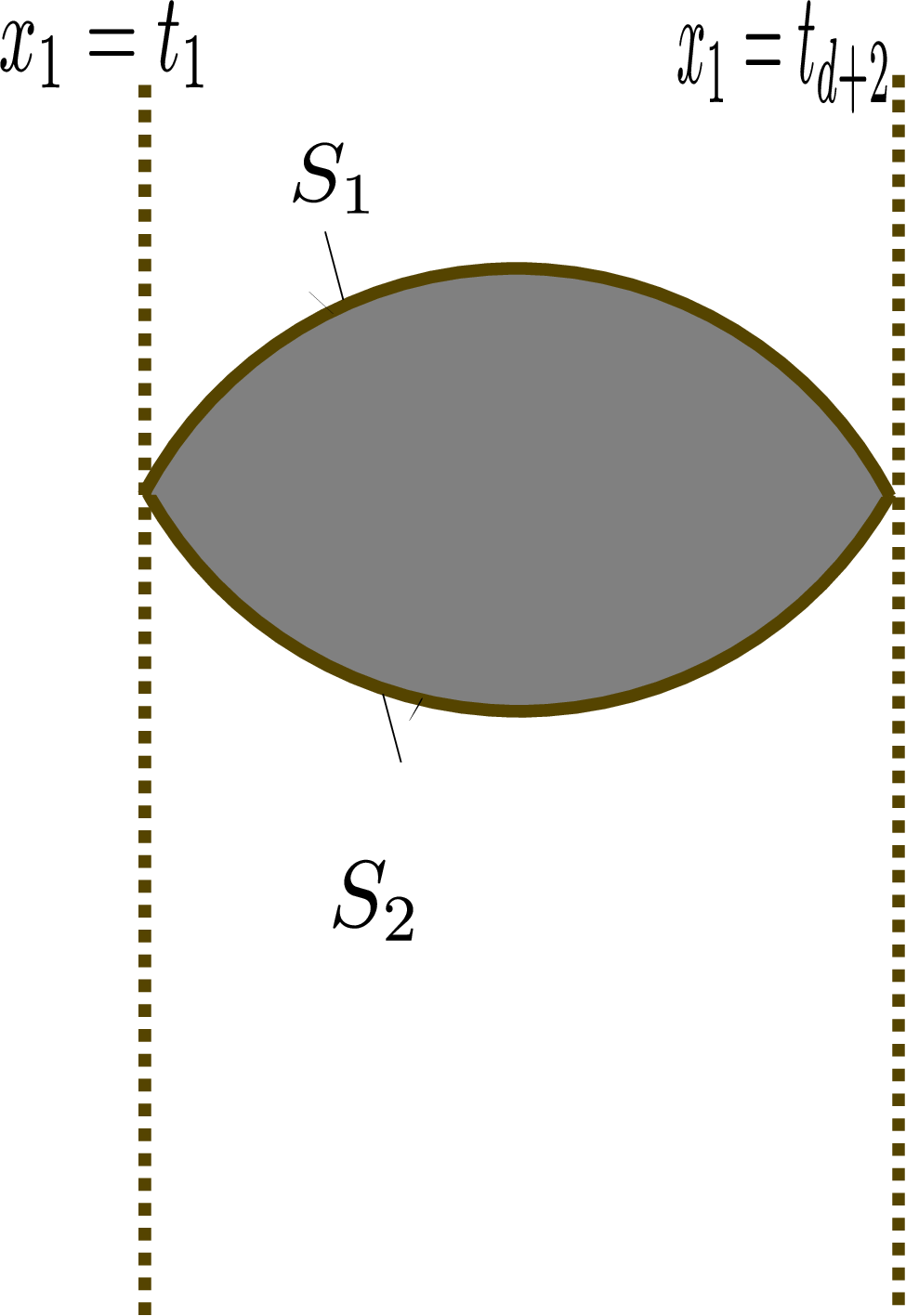}
	\caption{The first two cylinders $S_1$ and $S_2$ of $C_R$, colored in black, and the region $D_{S_1} \bigcap D_{S_2}$ with $\partial D_{S_i}=S_i$, colored in gray, in ${\mathbb{R}}^2=\{(x_1,x_2) \mid (x_1,x_2) \in {\mathbb{R}}^2\}$, where $D_{S_i}$ is a cylinder of a disk $D_R$ with $\partial D_R=C_R$.}
	\label{fig:1}
\end{figure}
\item The remaining cylinders $S_j$ ($j \geq 3$) of the circle are as presented in Figure \ref{fig:2}.  
These are defined for each depth $1 \leq i \leq d$ for vertices of the rooted tree. For each $i$, $\{(x_1,x_{i+2}) \mid (x_1,x_{i+2}) \in {\mathbb{R}}^2\}$ is depicted, the number of the chosen cylinders is $n_{c,i}+1$, and the depicted cylinders of the circles are in the plane $\{(x_1,x_{i+2}) \mid  (x_1,x_{i+2}) \in {\mathbb{R}}^2\}$, circles whose radii are $R$.
We can have a cylinder $D_{S_j}$ of the disk $D_R$ with $\partial D_{S_j}=S_j$ for each cylinder $S_j$ ($j \geq 3$) of the circle $C_R$. In this situation, a cylinder $E_{S_j}$ of the complementary set of the disk $D_R$ with $\partial E_{S_j}=S_j$ is chosen instead of the cylinder of the disk $D_R$ with $\partial D_{S_j}=S_j$. 
The region represented as the intersection ${\bigcap}_{j=1}^{{\Sigma}_{i=1}^d (n_{c,i}+1)} E_{S_{j+2}}$ is represented implicitly in FIGURE \ref{fig:2}.
%$D_{S_1} \bigcap D_{S_2}$
\begin{figure}
	\includegraphics[width=80mm,height=80mm]{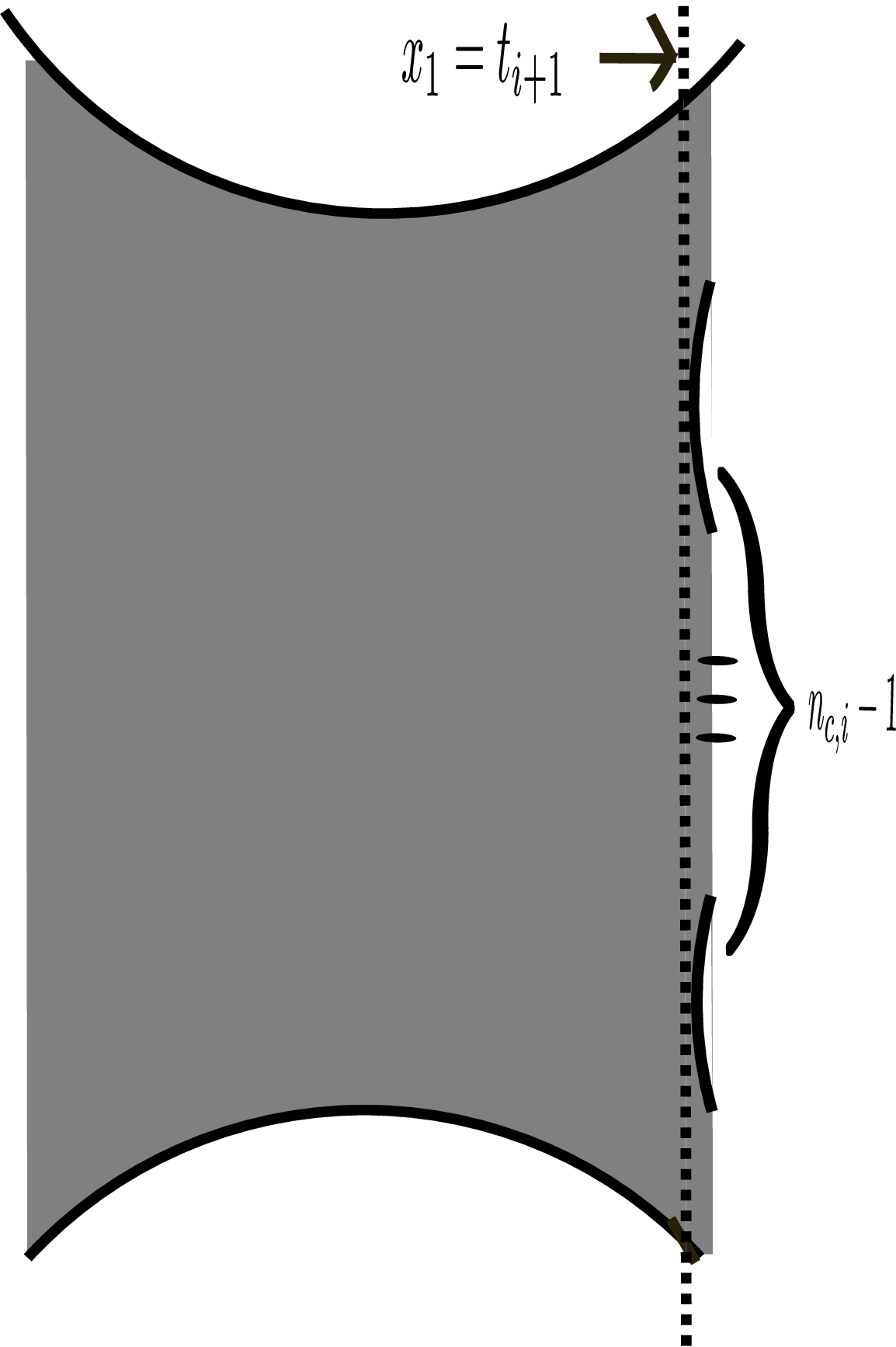}
	\caption{The remaining cylinders $S_j$ ($j \geq 3$) of the circle $C_R$ defined for each depth $1 \leq i \leq d$ for vertices of the rooted tree. The region represented as the intersection ${\bigcap}_{j=1}^{{\Sigma}_{i=1}^d (n_{c,i}+1)} E_{S_{j+2}}$ is represented in gray.}
	\label{fig:2}
\end{figure}

\end{itemize}

We can define a bounded, connected and open region $D$ surrounded by these cylinders as $D:=D_{S_1} \bigcap D_{S_2} \bigcap {\bigcap}_{j=1}^{{\Sigma}_{i=1}^{d} (n_{c,i}+1)} E_{S_{j+2}}$. 
We check that the pair $(D,\{S_j\}_{j=1}^{{\Sigma}_{i=1}^{d} (n_{c,i}+1)+2})$ satisfies Definition \ref{def:2}. 
By an elementary argument, the condition (\ref{def:2.1}) is satisfied. We check the condition (\ref{def:2.2}). For this, as in \cite{kitazawa5, kitazawa6, kitazawa7}, for each point $p$ in
${\overline{D}}^{{\mathbb{R}}^k}-D$, we consider a normal vector of each hypersurface $S_j$ containing $p$ at $p$ in ${\mathbb{R}}^k$ and check that these normal vectors are mutually independent.
Hereafter, let $e_j$ denote the $k$-dimensional real vector the value of whose $j$-th component is $1$ and the values of whose remaining components are $0$. The family of normal vectors of hypersurfaces $S_j$ at the point $p \in S_j$ can be chosen as follows.
\begin{itemize}
\item If $p$ is in either $S_1$ or $S_2$, then the family contains exactly one vector of the form $t_1e_1+t_2e_2$ with $t_1 \neq 0$ and $t_2 \neq 0$.  If $p$ is in both $S_1$ and $S_2$, then the family contains one vector of the form $t_1e_1+t_{2,1}e_2$ and exactly one vector of the form $t_1e_1+t_{2,2}e_2$ with $t_{2,1}>0$, $t_{2,2}<0$ and $t_1>0$. If $p$ satisfies neither of these two, the family does not contain a normal vector of such a form. 
\item If $p$ is in some $S_{i_j}$ chosen for the depth $1 \leq i \leq d$ for vertices of the rooted tree, then $p$ is in exactly one of such a cylinder $S_{i_j}$ of the circle $C_R$ by our choice and exactly one vector of the form $t_1e_1+t_ie_{i+2}$ is contained in the family with $(t_1,t_i) \neq (0,0)$. If $p$ is not contained in any $S_{i_j}$ chosen for the depth $1 \leq i \leq d$ for vertices of the rooted tree, such a vector is not contained. Furthermore, by our choice, at most one of vectors of the form here is of the form $t_1e_1$.
\end{itemize}
We can see that the family of the normal vectors of $S_j$ at $p$ are mutually distinct.

From the construction, this completes the proof.
\end{proof}
Main Theorem \ref{mthm:2} is a kind of additional theorems. 
\begin{MainThm}
\label{mthm:2}
Let $d_1>0$ and $d_2>0$ be integers. 
Let $\{t_j\}_{j=1}^{d_1+d_2+1}$ be an increasing sequence of real numbers{\rm :} $t_j<t_{j+1}$ for $1 \leq j \leq d_1+d_2+1$.
Let $\{n_{c,1,i}\}_{i=1}^{d_1}$ and $\{n_{c,2,i}\}_{i=1}^{d_2}$ be sequences of integers greater than $1$.
Let
$\overrightarrow{G}$ be a digraph obtained in the following way.
\begin{itemize}
\item Choose a balanced tree of type $(\{n_{c,1,i}\}_{i=1}^{d_1},d_1)$, let $v_{0,1}$ denote its root, and reverse the orientation.
\item Choose another balanced tree of type $(\{n_{c,2,i}\}_{i=1}^{d_2},d_2)$ and identify $v_{0,1}$ and the root $v_{0,2}$ of this balanced tree of type $(\{n_{c,2,i}\}_{i=1}^{d_2},d_2)$.
\end{itemize}
Then an RA-region $(D,\{S_j\}_{j=1}^{l}):=(D,\{S_j\}_{j=1}^{{\Sigma}_{i=1}^{d_1} (n_{c,1,i}+1)+{\Sigma}_{i=1}^{d_2} (n_{c,2,i}+1)+2})$ surrounded by cylinders $S_j$ of a circle $C_R:=\{x=(x_1,x_2) \mid ||x||=R\}$ in ${\mathbb{R}}^2$ whose radius is any real number $R$ greater than a sufficiently large positive number $R_0$ can be chosen in such a way that the Poincar\'e-Reeb digraph $\overrightarrow{{{PR}_{(D,\{S_j\}_{j=1}^l)}}_{\overline{{\pi}_{(D,\{S_j\}_{j=1}^l)}}}}$ is isomorphic to $\overrightarrow{G}$ and that the image of its vertex set by the map $\overline{{\pi}_{(D,\{S_j\}_{j=1}^l)}}:{PR}_{(D,\{S_j\}_{j=1}^l)} \rightarrow \mathbb{R}$ is $\{t_j\}_{j=1}^{d_1+d_2+1}$.
\end{MainThm}
\begin{proof}
We choose as in Figures \ref{fig:3}, \ref{fig:4} and \ref{fig:5} instead in the case of Main Theorem \ref{mthm:1} and its proof. We can prove similarly.
\begin{figure}
	\includegraphics[width=80mm,height=80mm]{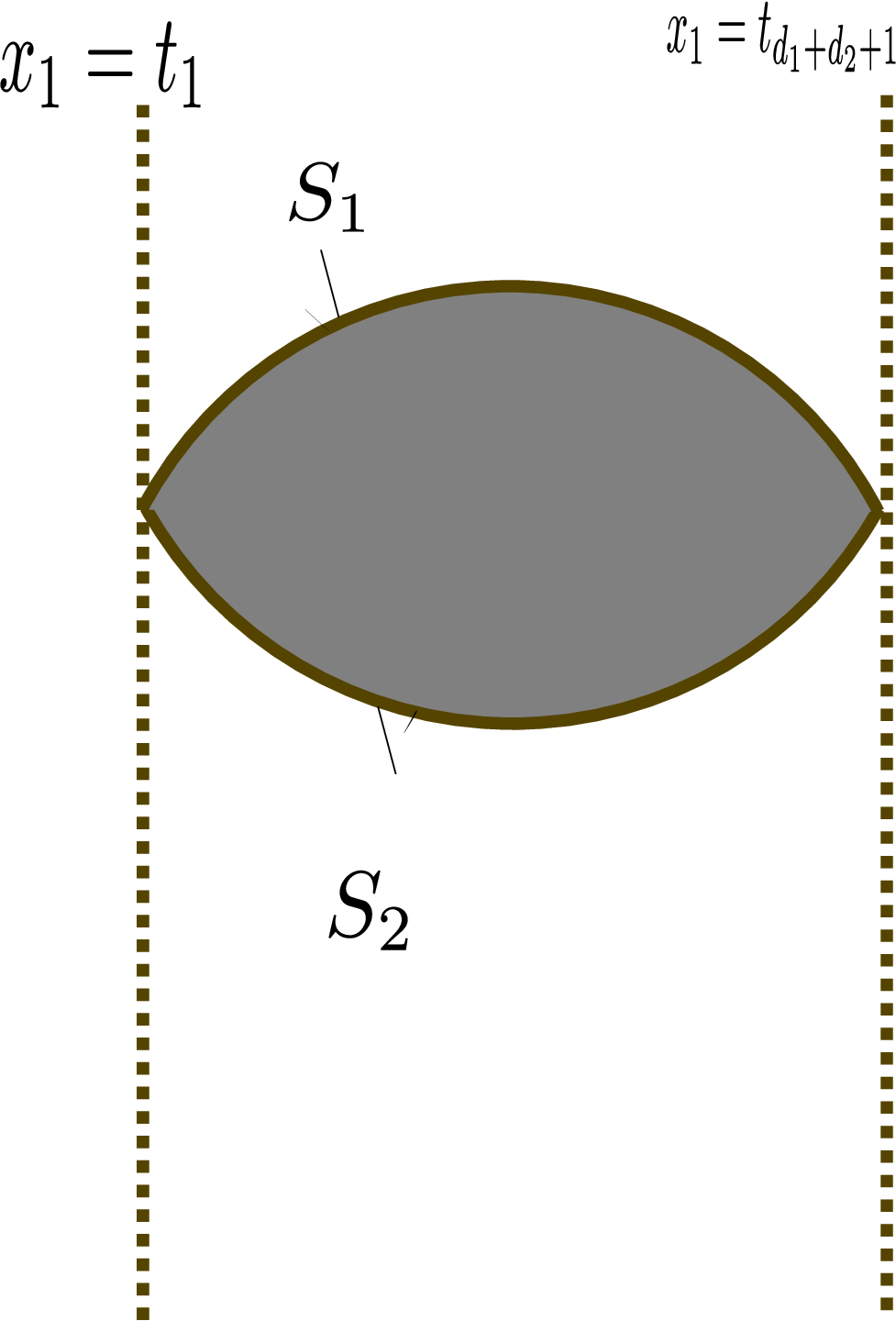}
	\caption{The cylinders $S_1$ and $S_2$ of $C_R$.}
	\label{fig:3}
\end{figure}
\begin{figure}
	\includegraphics[width=80mm,height=80mm]{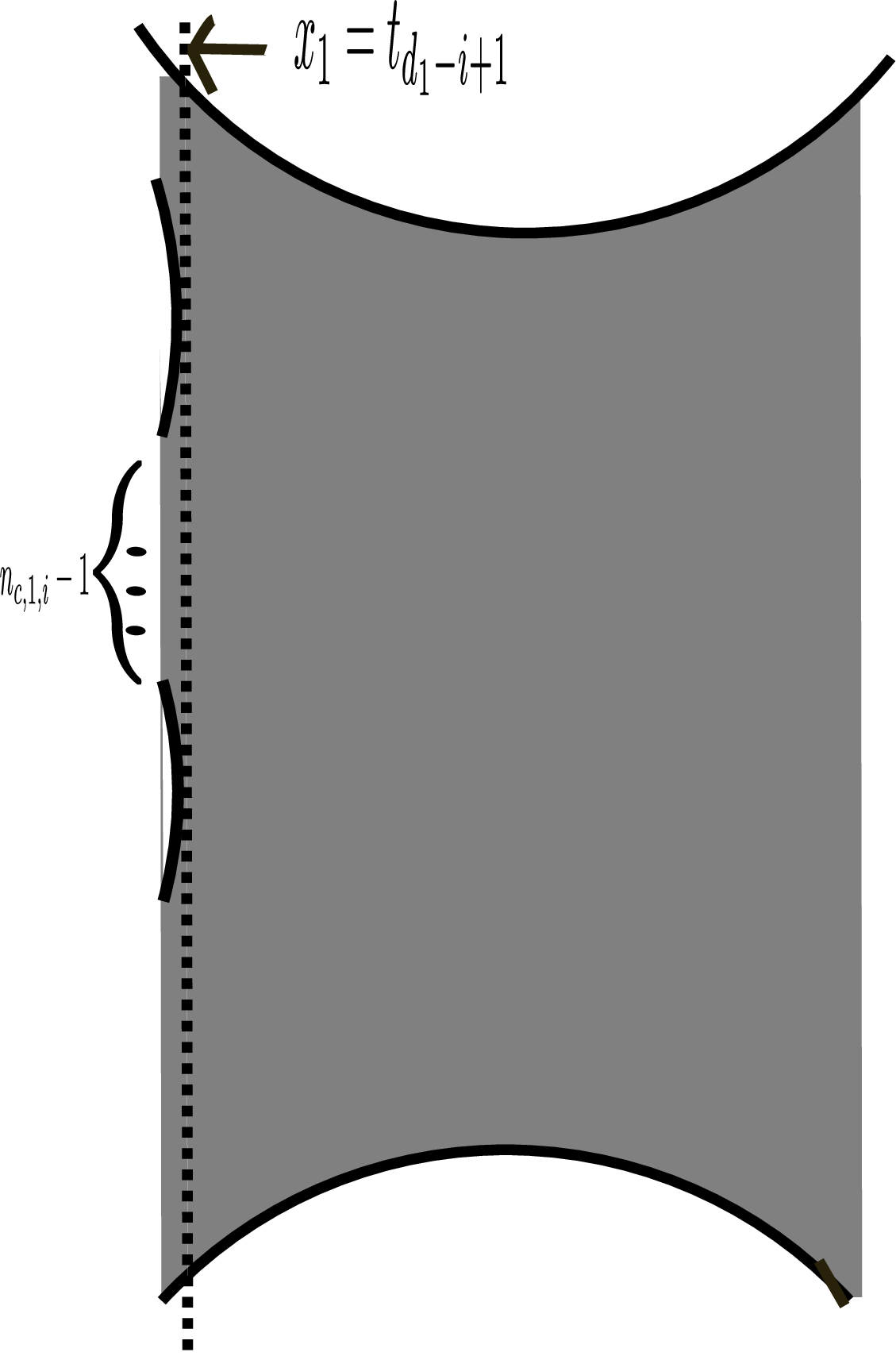}
	\caption{The remaining cylinders $S_j$ ($3 \leq j \leq 2+{\Sigma}_{i=1}^{d_1} (n_{c,1,i}+1)$) of the circle $C_R$ defined for the balanced tree of type $(\{n_{c,1,i}\}_{i=1}^{d_1},d_1)$ with the orientation reversed.}
	\label{fig:4}
\end{figure}
\begin{figure}
	\includegraphics[width=80mm,height=80mm]{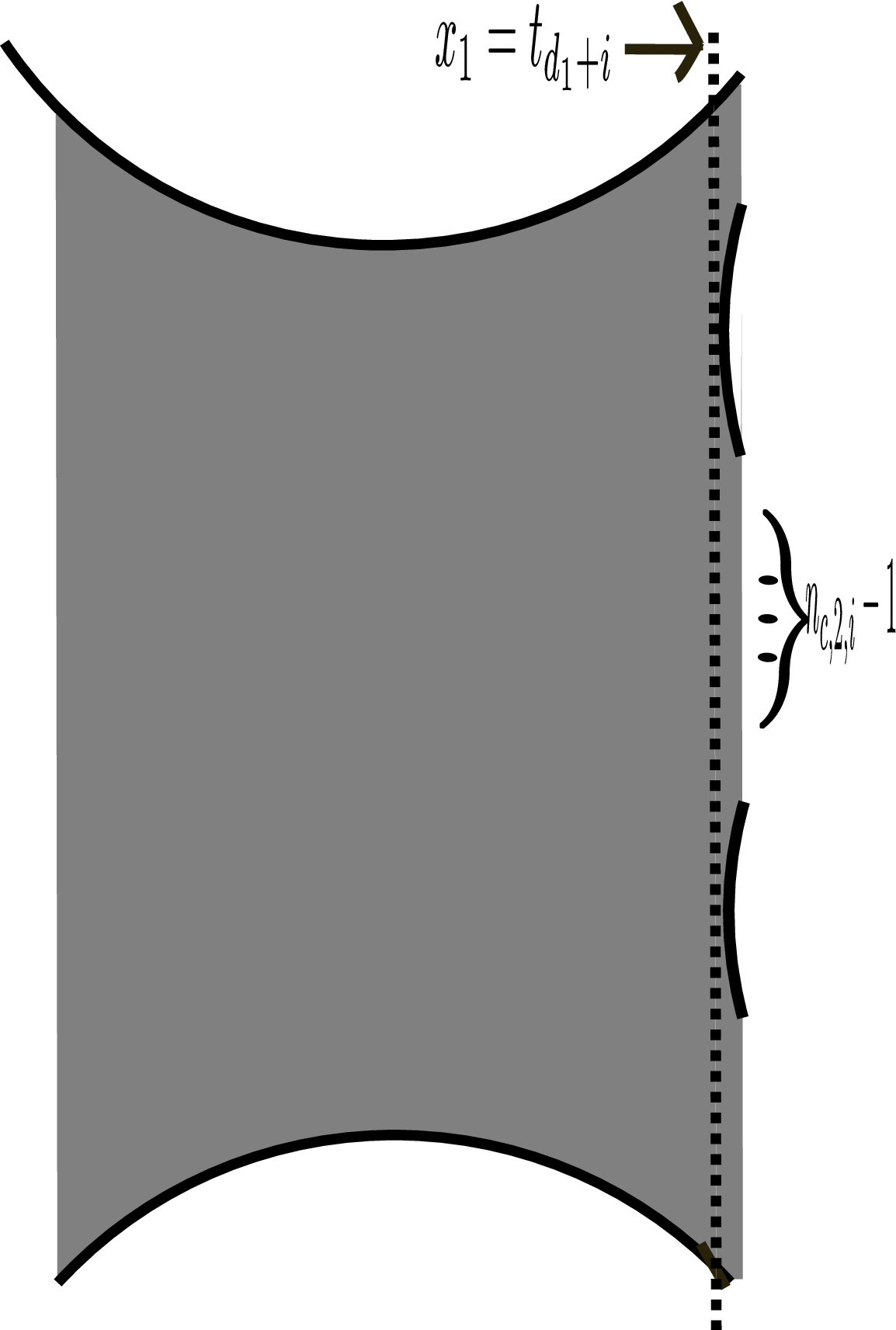}
	\caption{The remaining cylinders $S_j$ ($j \geq 3+{\Sigma}_{i=1}^{d_2} (n_{c,1,i}+1)$) of the circle $C_R$ defined for the balanced tree of type $(\{n_{c,2,i}\}_{i=1}^{d_2},d_2)$.}
	\label{fig:5}
\end{figure}
\end{proof}
\begin{Rem}
We drop the condition $n_{c,1,1}>1$ on the balanced tree of type $(\{n_{c,1,i}\}_{i=1}^{d_1},d_1)$ and let $d_1=1$, $d_2=d$ and $n_{c,1,1}=1$. This is Main Theorem \ref{mthm:1}.
\end{Rem}
% \cite{masumotosaeki} generalizes the pioneering work of \cite{sharko}. \cite{sharko} constructs nice smooth functions on closed surfaces and \cite{masumotosaeki} extends this to arbitrary finite graphs.
%Later, for example, \cite{martinezalfaromezasarmientooliveira, michalak} have set explicit problems and solved. Before the study \cite{kitazawa1} of the author, functions are, ones on closed surfaces or ones preimages containing no singular points of which are disjoint unions of spheres essentially.

%However we can apply some important cases. For example \cite{saeki2} considers very general cases and we cannot use analytic functions for construction there.
\section{Some remark on Main Theorems.}
\subsection{Some real algebraic maps.}
In Definition \ref{def:2}, let $d(j)>0$ be a positive integer defined for each integer $1 \leq j \leq l$. We can define the set $M_{(D,\{f_j\}_{j=1}^l)}=\{(x,y)=(x,\{y_j\}_{j=1}^{l} \in {\mathbb{R}}^k \times {\prod}_{j=1}^{l} {\mathbb{R}}^{d(j)+1}= {\mathbb{R}}^{l+k+{\Sigma}_{i=1}^{l} d(j)} \mid f_j(x)-{||y_j||}^2=0, 1 \leq j \leq l\}$. According to \cite[Main Theorem 1]{kitazawa3}, this is regarded as a connected component of the zero set of the real polynomial map (we can define naturally) and a real algebraic manifold. 
By restricting ${\pi}_{l+k+{\Sigma}_{j=1}^{l} d(j),k}$ to $M_{(D,\{f_j\}_{j=1}^l)}$, a real algebraic map is defined. From the viewpoint of singularity theory, this map is, around each of its singular point, locally like a so-called {\it moment map}. Moment maps are fundamental and important in geometry of so-called toric symplectic manifolds and see the textbook \cite{buchstaberpanov} and see also \cite{delzant}.  

We can also have the function ${\pi}_{l+k+{\Sigma}_{i=1}^{l} d(i),1}$ and this is shown to be a so-called {\it Morse-Bott} function in \cite{kitazawa3}. Morse-Bott functions are generalizations of Morse functions: see the textbook \cite{banyagahurtubise} and see also \cite{bott}.

We can also have the fact that the Reeb digraph $\overrightarrow{{R_f}_{\bar{f}}}$ and the Poincar\'e-Reeb digraph of $(D,\{S_j\}_{j=1}^l)$ are isomorphic.  

Arguments here can be generalized to some extent as in the original preprint. We introduce in a specific way for simplicity.

\subsection{Some cases related to our main result.}
In \cite{kitazawa6}, we have obtained an explicit variant of Main Theorems for cylinders of a connected component $C=\{(x_1,x_2) \in {\mathbb{R}}^2 \mid x_1x_2-1=0, x_2>0 {\rm  (}x_2<0{\rm )}\}$ of a {\it hyperbola} of the form $\{(x_1,x_2) \in {\mathbb{R}}^2 \mid x_1x_2-1=0\}$. 
This is to have a real algebraic function as presented in the previous subsection in such a way that the Reeb digraph is a path digraph on $d+1$ vertices and that the function satisfies arbitrary prescribed conditions on (non-)compactness of level sets for points in the interiors of each edge. This result follows the first answer to this explicit case, presented in \cite{kitazawa5}. For this, see also \cite{kitazawa2}, a related study in the differentiable (smooth) situation motivating the author to study the real algebraic version. See \cite{kitazawa14} and see also \cite{kitazawa13} for previous studies on construction of Reeb digraphs or Poincar\'e-Reeb digraphs being also trees. 

Last, for cases where $S_j$ are cylinders of circles, see also \cite{kitazawa7}. This gives examples of Reeb graphs of real algebraic functions which are not planar.

\end{document}